\pdfoutput=1
\documentclass{article}
\pdfoutput=1

\usepackage{url}
\usepackage{lineno}
\usepackage{enumitem}
\usepackage{times,amsfonts,amsmath,mathtools,epsfig,fixltx2e,amsthm,amssymb}
\usepackage{graphics,graphicx,epstopdf}
\usepackage{caption}
\usepackage{cite}
\usepackage{multirow}

\newcommand{\reals}{{\mathrm{I\kern-.2em R}}}
\newcommand{\complex}{{\mathrm{C\kern-.6em C}}}
\newcommand{\field}{{\mathrm{I\kern-.2em F}}}
\newcommand{\expectation}{{\mathrm{I\kern-.2em E}}}

\newcommand{\calD}{{\cal D}}

\newcommand{\calH}{{\cal H}}

\newcommand{\calP}{{\cal P}}

\newcommand{\m}{m^\prime}

\newcommand{\vx }{{\bf x}}

\newcommand{\qedwhite}{\hfill \ensuremath{\Box}}

\newcommand{\basispln}{I_{p;l,n}}

\newcommand{\basisplm}{I_{p;l,m}}
\newcommand{\scalarbasispln}{I_{p;l,n}}

\newcommand{\Nirrep}{N_{\rm rep}}
\newcommand{\Ngroup}{N_g}
\newcommand{\Npl}{N_{p;l}}

\makeatletter
\DeclareRobustCommand*\cal{\@fontswitch\relax\mathcal}
\makeatother
\title{Real basis functions of polyhedral groups}

\author{
	Nan Xu\thanks{School of Electrical and Computer Engineering, Cornell University, Ithaca, NY, USA, ID ({nx25@cornell.edu}).}
}

\date{March 16, 2016}

\begin{document}

\maketitle
\begin{abstract}
The basis of the identity representation of a polyhedral group is able to describe functions with symmetries of a platonic solid, i.e., 3-D objects which geometrically obey the cubic symmetries. However, to describe the dynamic of assembles of heterogeneous 3-D structures, a situation that each object lacks the symmetries but obeys the symmetries on a level of statistics, the basis of all representations of a group is required. While those 3-D objects are often transformed to real functions on $L_2$ space, it is desirable to generate a complete basis on real space. This paper deduces the existence of a basis on real space for each polyhedral group, and introduces a novel approach to explicitly compute these real basis functions, of which properties are further explored.
\end{abstract}

\section{Introduction}
\label{sec:intro}
Basis functions of identity representation (rep) of a polyhedral group are able to describe any function adapted with symmetry of a platonic solid, which includes Tetrahedral, Octahedral and Icosahedral symmetry, whereas the more general basis functions of all reps, which have certain transformation property maintained, are able to describe any function on the $L^2$ space. A popular problem in x-ray crystallography and biophysics is to describe geometric behavior of a class of virus, namely ``spherical'' virus, which geometrically exhibits the symmetry of a polyhedral group, i.e., Tetrahedral symmetry, Octahedral symmetry or Icosahedral symmetry. The angular behavior of these bio-nanomachines can be expressed an orthonormal expansion of basis functions of the corresponding polyhedral group. A recent interest in structural virology is to describe the geometric behavior of an assemble of virus particles, in which each individual lack the symmetry but only obey the symmetry on the first and second order statistics. This is a sophisticated but realistic scenario, which is able to accommodate the heterogeneity of virus particles during their maturation process. To describe such sophisticated scenario, the explicit computation of a complete orthonormal set (C.O.N.S) of basis functions of a polyhedral group that spans on the $L^2$ space is required. Due to the fact of all image data is real-valued, it is desirable to compute real basis functions of polyhedral groups.

There has been extensive study on basis functions of a polyhedral group~\cite{fox1970construction,muggli1972cubic,AltmannPCambPhilSoc1957,MeyerCandJMath1954,CohanPCambPhilSoc1958,ElcoroPerezMatoMadariaga1994,HeuserHofmannZNaturforsch1985,JackHarrison1975,KaraKurkiSuonioActaCryst1981,LaporteZNaturforschg1948,conte1984determination,raynal1984determinationIII,ZhengDoerschukActaCryst1995,FernandoJComputPhys1994}. Some focus on basis functions of the identity representation to describe specific symmetric functions, and others focus on the more general problem of basis functions of all reps of a point group~\cite{CohanPCambPhilSoc1958,conte1984determination,raynal1984determinationIII}. In the majority of existing literature, basis functions of a symmetry group have been generated as a linear combination of spherical harmonics of a single degree~\cite{AltmannPCambPhilSoc1957,Altmann1962,Altmann1965,MuellerPriestley1966,Puff1970,muggli1972cubic,fox1970construction,ZhengDoerschukActaCryst1995,zheng2000explicit}, because of the importance of rotations and the relative simplicity of rotating spherical harmonics. This concept has been widely applied for  the last four decades, e.g., the fast rotation function~\cite{CrowtherFastRotation1972}. Other work express the basis functions of a polyhedral group as multipole expansions in the rectangular coordinates~\cite{Jahn469,HECHT1961}. 
Previous work uses a variety of techniques and often has a restriction on the value $l$ of the spherical harmonics~\cite{AltmannPCambPhilSoc1957,Altmann1962,Altmann1965,MuellerPriestley1966,Puff1970,CohanPCambPhilSoc1958,muggli1972cubic,conte1984determination,raynal1984determinationIII}. For instance, Refs.~\cite{AltmannPCambPhilSoc1957,Altmann1962,Altmann1965} consider a range of point groups and use the techniques of projection operators and Wigner $D$ transformations to compute basis functions up to degree $l=12$, where Ref.~\cite{CohanPCambPhilSoc1958} uses similar techniques restricted to the icosahedral group to provide basis functions up to degree $l=15$.
Refs.~\cite{MuellerPriestley1966,Puff1970} use the method of representation transformation to compute the invariant basis functions of the cubic group, which can compute up to degree $l=30$; the work of \cite{muggli1972cubic} extends this computation to all irreps basis functions. 
Refs. \cite{conte1984determination,raynal1984determinationIII} proposed a method to deriving all irreps basis functions of the cubic and the icosahedral groups for a specific degree $l$. However, for computation which needs all irreps basis functions of a large range of $l$'s (e.g., from 1 to 55), the one by one derivation is cumbersome. Later work~\cite{fox1970construction,ZhengDoerschukActaCryst1995,zheng2000explicit} release this degree restriction and allow for the computation of the invariant basis function of any polyhedral groups. We note that the recursions in~\cite{ZhengDoerschukActaCryst1995,zheng2000explicit} appear to be unstable in computation experiments. 
\par 
In our paper, we derive the explicit form for the complete basis functions of all reps of the tetrahedral, the octahedral and the icosahedral groups. To compute the real orthogonal basis, we seek real projection operators and then apply them to a real spherical harmonics. Real-valued projection operators require real irreducible representation matrices. The standard approach, e.g., Young diagrams~\cite{fulton1997young}, are able to generate irreducible representation matrices. Benefit from the existing complex irreducible representation matrices\cite{Aroyo2006b,LiuPingChenJMathPhys1990}, we establish a computational method based on the Frobenious-Schur indicator~\cite{GanevSchurIndicator2011} to transform these complex matrices to real irreducible representation matrices. The solution depends on a matrix eigenvalue problem of dimension equal to twice the dimension of the irreducible representation which is practical for numerical computation for a finite irreducible symmetry group. We also demonstrate that the absence of a real irreducible representation implies that real basis functions do not exist.  
\par 
For particles with exactly the same cubic symmetry, its angular behavior can be described by the basis functions of identity representation of that polyhedral group \cite{zheng2000explicit}. Specifically, basis functions of identity representation of Icosahedral group has been computed as linear combinations of spherical harmonics~\cite{zheng2000explicit}, which have been used to study the structure of virus particles with icosahedral symmetry recently~\cite{YiliZhengQiuWangDoerschukJOSA2012,QiuWangMatsuiDomitrovicYiliZhengDoerschukJohnsonJSB2012}.
The approach in \cite{zheng2000explicit} can be generalized to compute angular basis functions for particles with non-icosahedral symmetry, i.e., tetrahedral and octahedral symmetries, which have not been widely explored. However, during the maturation process of a virus, particles have continuous variability resulting in heterogeneous structures without obeying any symmetry. Therefore, it is more realistic to assume that each particle lack the specific symmetry during the life-cycle but obey the symmetry on a level of statistics. Then, basis functions of all representations of a polyhedral group is required to describe such a sophisticated situation, which is the motivation of this paper.
\par
The sections of the paper is arranged in the following. Section~\ref{sec:realirrepmatrices} computes the real irreducible representation matrices from the complex ones. Section~\ref{sec:computebasisfunctions} introduces the computation of the real basis functions and their properties. Section~\ref{sec:results} shows the results of real basis functions.

\section{Real irreducible representations of a polyhedral group}\label{sec:realirrepmatrices}
Tetrahedral group $T$, Octahedral group $O$ and Icosahedral group $I$ are the three symmetry groups of the platonic solids\footnote{Among platonic solids, cube and octahedron are dual polyhedron of each other, so as dodecahedron and icosahedron. The former two obey octahedral symmetries, and the later two obey icosahedral symmetries.}, which are referred as the polyhedral groups. Their irreducible representations (irreps, ``rep'' for ``representation'') on complex space are well studied~\cite{Aroyo2006b,LiuPingChenJMathPhys1990}. In this section, an approach to transforming the possible irreps from the complex space to the real space is introduced. 

A polyhedral group $G$ of order $N_g$ has $N_{rep}$ irreps. The $p$th irrep, for $p=1,2 ..., N_{rep}$, can be represented by a set of $d_p\times d_p$ unitary matrices, denoted by $\Gamma_c^{p}\in\mathbb{C}^{d_p\times d_p\times N_g}$, in which each matrix $\Gamma_c^{p}(g)\in\mathbb{C}^{d_p\times d_p}$ represents a group operation $g\in G$.
The corresponding values of $N_{g}$, $N_{rep}$ and $\{d_1, ..., d_{N_{rep}}\}$ for each symmetry group of platonic solids are tabulated in Table \ref{tab:symGp}. 
\begin{table}[h]
	\begin{center}
		\begin{tabular}{|c|c|c|c|}
			\hline
			Symmetry Groups & $N_g$ & $N_{rep}$ & $d_p$  \\
			\hline
			Tetrahedral & 12 & 4 & \{1, 1, 1, 3\} \\
			\hline
			Octahedral & 24 & 5 & \{1, 1, 2, 3, 3\} \\
			\hline
			Icosahedral & 60 & 5 & \{1, 1, 2, 3, 3\}\\
			\hline
		\end{tabular}
	\end{center}
	\caption{\label{tab:symGp} Properties of the three polyhedral groups.}
\end{table}

Complex-valued unitary irrep matrices for Tetrahedral group $I$ and Octahedral group $O$ are available on the Bilbao Crystallographic Server~\cite{Aroyo2006b}. Complex-valued unitary irrep matrices for Icosahedral group is provided by Liu, Ping, and Chen~\cite{LiuPingChenJMathPhys1990}.

\subsection{Potentially real irreps}
Let $\chi^{p}(g)$ be the character of element $g\in G$ in the $p$th rep, i.e., the trace of matrix $\Gamma_c^{p}(g)$. Based on~\cite[p.~129, Theorem~III]{Cornwell1984}, the $p$th irrep $\Gamma_c^p$ is {\it potentially real}, if 
\begin{equation}
1/N_g\sum_{g\in G}\chi^p(g)^2=1/N_g\sum_{g\in G}\chi^p(g^2)=1. 
\label{eq:RealIrredRepCond}
\end{equation}

That is, if Eq.~\ref{eq:RealIrredRepCond} holds, there existing real-valued $p$th irrep matrices, denoted by $\Gamma_r^p\in\mathbb{R}^{d_p\times d_p\times N_g}$, which is similar to the complex-valued irrep matrices $\Gamma_c^p$~\cite[p.~128, Theorem~II]{Cornwell1984}. 
Specifically, there exists some unitary transformation matrix $S^p\in\mathbb{C}^{d_p\times d_p}$ for $p\in\{1,2,...,5\}$, such that, 
\begin{equation}
\Gamma_r^p(g)=(S^p)^H\Gamma_c^p(g)S^p\label{eq:SimilarTransformation}\text{~and $\Gamma_r^p(g)\in\mathbb{R}^{d_p\times d_p}$ for all $g\in G$.}
\end{equation}  Note that, for each $p$ value, the obtained $\Gamma_r^p$ has the same multiplication table as $\Gamma_c^p$ does, of which both are homomorphic to group operations. 
By testing Eq.~\ref{eq:realProperty} on all polyhedral groups, we have Eq. \ref{eq:RealIrredRepCond} hold for all except for the 2nd and the 3rd irreps of Tetrahedral group, which have the left hand side of Eq. \ref{eq:RealIrredRepCond} equal to zero. Following Theorem~III in~\cite[p.~129]{Cornwell1984}, we have
\begin{quote}
	Corollary 1: {\sl For Octahedral and Icosahedral groups, all irreps are potentially real. For Tetrahedral group, only the 1st and the 4th irreps are potentially real, but the 2nd and the 3rd irreps can only be complex.
	}
\end{quote}


\subsection{Computation of real irrep matrices from complex irrep matrices}\label{sec:FindTransformSp}
If we have complex irrep matrices that are potentially real, one then is able to compute real irrep matrices by Eq.~\ref{eq:SimilarTransformation}. The property of real-valued matrices is maintained by taking the complex conjugate of Eq. \ref{eq:SimilarTransformation} and equaling it to the original, i.e. $\overline{\Gamma_r^p(g)}=\Gamma_r^p(g)$ or
\begin{alignat}{2}\label{eq:realProperty}
\overline{(S^p)^H\Gamma_c^p(g)S^p}&=(S^p)^H\Gamma_c^p(g)S^p,\text{~so that~}\nonumber\\
\overline{\Gamma_c^p(g)}=(S^p{S^p}^{T})^{H}\Gamma_c^p(g)(S^p{S^p}^{T})&=(S^p{S^p}^{T})^{-1}\Gamma_c^p(g)(S^p{S^p}^{T}) \text{~for all~} g\in G.
\end{alignat}
The second equal sign in the above equation is due to the unitarity of matrix $S^p$.
In other words, the complex conjugate $\overline{\Gamma_c^p}$ becomes similar to the original ${\Gamma_c^p}$ due to its similarity to the real irrep $\Gamma_r^p$~\cite[p.~127]{Cornwell1984}. Followed by Theorem~II in~\cite[p.~128]{Cornwell1984}, a non-singular matrix $Z^p$, which satisfies the form
\begin{equation}\overline{\Gamma_c^p}=(Z^p)^{-1}\Gamma_c^p(g)Z^p \text{for all $g\in G$, where }\overline{Z^p}Z^p=c_z^pI, \label{eq:complexSimilarity}\end{equation}for $c_z^p\in\mathbb{R}^+$ and an $d_p\times d_p$ identity matrix $I$, can be determined by \begin{equation} Z^p=\frac{1}{N_g}\sum_{g\in G} \Gamma^p_c(g)A^p(\overline{\Gamma^p_c(g)})^{-1} =\frac{1}{N_g}\sum_{g\in G} \Gamma^p_c(g)A^p{\Gamma^p_c(g)}^{T}\label{eq:Zconstruction}\end{equation} 
for some matrix $A^p\in \mathbb{C}^{d_p\times d_p}$.
Consider the matrix $C^p=Z^p/\sqrt{c_z^p}=S^p{S^p}^{T}$. Then, $C^p\in\mathbb{C}^{d_p\times d_p}$ is unitary.
A simple example of $C^p$ is to consider it as symmetrical (but not hermitian), i.e., ${C^p}={C^p}^T$, which can be computed from a random symmetrical matrix $A^p\in\mathbb{C}^{d_p\times d_p}$ by Eq. \ref{eq:Zconstruction}.  In such case, $C^p$ satisfies the relation of
\begin{equation}\overline{\Gamma_c^p}=\overline{C^p}\Gamma_c^p(g)C^p \text{ for all $g\in G$, where }\overline{C^p}C^p=I.\label{eq:complexSymSimilarity}\end{equation} 
Takagi Factorization (and Corollary 4.4.6 in Ref.~\cite[p.~207]{HornJohnson1985}) guarantees the existence of such unitary $S^p$. One approach to computing $S^p$ is introduced in next paragraph. 

Consider the real and complex form of $C^p$, i.e., $C^p=C^p_R+iC^p_I$ with $C^p_R, C^p_I\in\mathbb{R}^{n\times n}$. Then, it has the real representation matrix $B^p=\begin{bmatrix}
C^p_R & C^p_I\\
C^p_I & -C^p_R
\end{bmatrix}\in\mathbb{R}^{2n\times 2n}$. Followed by Proposition 4.6.6 in Ref.~\cite[p.~246]{HornJohnson1985} and Lemma 1 in Appendix A, eigenvalues of $B^p$ appear in $(+1,-1)$ pairs. Let $\begin{bmatrix} X^p \\ -Y^p \end{bmatrix}$ be the matrix columned by the orthonormal eigenvectors of $B^p$ that associates to $+1$ eigenvalues. Then, $S^p$ can be chosen to be $S^p=X^p-iY^p$. Such $S^p$ has been verified with its underlying theory introduced in Appendix A. The software and numerical results are available upon requests.
\section{Real orthonormal basis functions of polyhedral groups}\label{sec:computebasisfunctions}
Let the $d_p$-dimensional vector function $I_{p;\zeta}(\vx/x)$ for $\zeta\in\{1,\dots,N_\zeta\}$ be an orthonormal basis function for the $p$th irrep of a polyhedral group. It is defined by a specific rotational operation from the group~\cite[p.~20]{Cornwell1984}, in
particular,
\begin{equation}
P(g)I_{p;\zeta}(\vx)=I_{p;\zeta}(R_g^{-1} \vx)=(\Gamma^p(g))^T
I_{p;\zeta}(\vx) \text{~for very $\zeta\in\{1,\dots,N_\zeta\}$},
\label{eq:BasisDefn}
\end{equation}
where $P(g)$ is the rotational operator representing the group operation $g\in G$, $R_g\in\mathbb{R}^{3\times 3}$ is the $g$th spatial rotation matrix for the symmetry group, and $\Gamma^p(g)$ is a $d_p\times d_p$ matrix representing the group operation $g$ in the $p$th irrep ($p\in\{1,\dots,\Nirrep\}$). From the definition, we can deduce the following (see Appendix B for the proof).
\begin{quote}
	Proposition 1: {\sl Real orthonormal basis functions, i.e., $I_{p;\zeta}\in\mathbb{R}^{d_p}$, generate real orthogonal irrep matrices of a group. }
\end{quote} 
Therefore, for the 2nd and 3rd irrep of Tetrahedral group, real-valued basis functions do not exist. In the reminder of this Section, we describe an approach to computing real orthonormal basis functions for Octahedral group, Icosahedral group, and the 1st and 4th irrep of Tetrahedral group.

A group theory approach to determine a basis function of a polyhedral group $G$ is to apply projection operators to a random functions of $L^2$ space~\cite[p.~93]{Cornwell1984}. 
Specifically, the projection operator
${\cal P}_{j,k}^p$ applied to a function $\psi(\vx)$ is a weighted sum
of rotational operators $P(g)$ applied to $\psi(\vx)$ where the weights are
matrix elements of an irrep and: 
\begin{equation}
{\cal P}_{j,k}^p
\psi(\vx)
=
\frac{d_p}{\Ngroup}
\sum_{g\in G}
\overline{\Gamma^p(g)_{j,k}}
P(g) \psi(\vx)
\label{eq:GroupProjectorOperation}
\end{equation}
where $\Gamma^p(g)_{j,k}$ for $j,k\in\{1,\dots,d_p\}$ is the $(j,k)$th element of the irrep matrix $\Gamma^p(g)$. One approach to determining real-valued polyhedral basis functions are to use real-valued irrep matrices, i.e., $\Gamma^p=\Gamma_r^p$, and the functions $\psi(\vx)$ with their rotational operators $P(g)$ that are real. 
\par
A basis of a point group can be constructed by projecting the point solid onto a sphere with the use the spherical harmonics~\cite{XuDoerschukICIP2015}. Spherical harmonics, denoted by $Y_{l,m}$ for $l\in\mathbb{N}$ and $m\in\{-l,\dots,l\}$, are sets of complex functions with simple rotational property,
$P(\cdot)Y^l
=
{D^l}^T Y^l
$, 
where $Y^l=[Y_{l,-l},\dots, Y_{l,l}]^T\in\mathbb{C}^{(2l+1)\times1}$ and $D^l\in\mathbb{C}^{(2l+1)\times (2l+1)}$ is the Wigner-D
matrix, in which the $(m,\m)$th entry is commonly referred as the Wigner D coefficient $D_{m,\m}^l$. Real spherical harmonics exist and form a complete orthonormal basis for $L^2$. To compute real basis function of a polyhedral group, it is natural to apply the projection operator to real spherical harmonics. 

With the valid real $P(\cdot)\psi(\vx)$ and valid real irrep matrices, according to Corollary 1, we can compute real basis functions for all irreps of Octahedral and Icosahedral groups, and for the 1st and 4th reps of Tetrahedral group. We have also shown the other direction of Proposition 1 for a polyhedral group.
\begin{quote}
	Proposition 2: {\sl Real basis functions of a polyhedral group exist if and only if the corresponding irrep is potentially real.}	
\end{quote}
\subsection{Rotation of real spherical harmonics}
\par
Real spherical harmonics and the rotation operation for real spherical harmonics can be transformed from the typical complex forms~\cite{collado1989rotation,Pinchon2007RealSpherical,aubert2013alternative}. Following Eq. 7 and 8 in~\cite{collado1989rotation}, real spherical harmonics, denoted by
$Z_{l,m}$ ($l\in\mathbb{N}$, $m\in\{-l,\dots, +l\}$), can be transformed by a $l\times l$ unitary matrix $U^l=$
\begin{eqnarray}
{\scriptsize
	\frac{1}{\sqrt{2}}
	\left[\begin{array}{ccccccccc}
	i &   &       &  &  &  & &   & 1\\
	& i &       &  &  &  & & 1 & \\
	&   &\ddots &  &  &  &\mathstrut^{.^{.^{.}}} & & \\
	&   &       &i &  &1 &  & & \\
	&   &       & &\sqrt{2}  & &  & & \\
	&   &       &i &  &-1 &  & & \\  
	&   &       -i& &  & &1  & & \\  
	&   i&       & &  & &  &-1 & \\  
	\mathstrut^{.^{.^{.}}}  &   &       & &  & &  & &\ddots\end{array} \right].
}
\end{eqnarray}

Let $Z^l=[Z_{l,-l}, \dots, Z_{l,l}]^T\in\mathbb{R}^{(2l+1)\times 1}$ be the vector of real spherical harmonics. 
Then, $Z^l={U^l}^T Y^l$ and the corresponding rotation matrix\footnote{There are two rotation matrix concepts here: $D^l\in\mathbb{C}^{(2l+1)\times(2l+1)}$ and $W^l\in\mathbb{C}^{(2l+1)\times(2l+1)}$ are the rotation matrices for complex and real spherical harmonics, respectively, whereas $R\in\mathbb{R}^{3\times 3}$ is the rotation matrix for the 3-D space coordinates, which is referred as symmetry rotation matrix in the later context.} is $W^l={U^l}^{-1} {D^l} {U^l}$. Hence, the rotational property becomes
\begin{equation}
P(\cdot)Z^l={W^l}^T Z^l=({U^l}^{-1} {D^l} {U^l})^T {U^l}^T Y^l=M^T Y^l
\label{Eq:realrotationMat}
\end{equation}
where $M^l={D^l} {U^l}$,
returns a $(2l+1)\times 1$ vector on real space. So, the rotation on real spherical harmonics can be computed as a linear transformation on complex spherical harmonics, which have been typically used. Choose $\psi(\vx)=Z_{l,m}(\vx/x)$ ($l\in\mathbb{N}$, $m\in\{-l,\dots, +l\}$). Then, 
\begin{eqnarray}
\label{eq:rotateYlmusingWignerD}
P(g)Z_{l,m}\left(\frac{\vx}{x}\right)
&=&Z_{l,m}\left(R_{g}^{-1} \frac{\vx}{x}\right)
=
\sum_{\m=-l}^{+l}
W_{\m,m}^l(g) Z_{l,\m}(\frac{\vx}{x})\\\nonumber
&=&
\sum_{\m=-l}^{+l}
M_{\m,m}^l(g) Y_{l,\m}(\frac{\vx}{x})
\end{eqnarray}
where $M_{\m,m}^l(g)=$
\begin{equation}\label{eq:NewRotationY}
\left\{
\begin{array}{lr}
\frac{i}{\sqrt{2}}\left(D_{\m,m}^l(g)+(-1)^{m+1} D_{\m,-m}^l(g) \right) & : m < 0\\
D_{\m,0}^l(g) & : m = 0\\
\frac{1}{\sqrt{2}}\left(D_{\m,-m}^l(g)+(-1)^m D_{\m,m}^l(g)\right) & : m > 0
\end{array}
\right..
\end{equation}

\subsection{Computation of Wigner-D coefficients}
\par
Computation of $M_{\m,m}^l$ (Eq.~\ref{eq:NewRotationY}) requires the knowledge of Wigner-D coefficients $D_{\m,m}^l(g)$ which depends on the Euler angles $(\alpha_g,\beta_g,\gamma_g)$ describing the rotation corresponding to each group operation $g$. These angles are not unique, since if the platonic solid is positioned in different orientations, the symmetry rotations (and therefore the Euler angles) are typically different. 

The symmetry rotation matrix $R_g$ is a function of the Euler angles $(\alpha_g,\beta_g,\gamma_g)$ in the following form~\cite[Eq. (4.43), p.65]{rosvall2008maps}:
\begin{equation}
{\scriptsize
	R_g=\begin{bmatrix}
	cos~\alpha_g cos~\beta_g cos~\gamma_g-sin~\alpha_g sin~\gamma_g & sin~\alpha_g cos~\beta_g cos~\gamma_g+cos~\alpha_g sin~\gamma_g & -sin~\beta_g cos~\gamma_g \\
	-cos~\alpha_g cos~\beta_g sin~\gamma_g-sin~\alpha_g cos~\gamma_g & -sin~\alpha_g cos~\beta_g cos~\gamma_g+cos~\alpha_g sin~\gamma_g & sin~\beta_g sin~\gamma_g \\
	cos~\alpha_g sin~\beta_g & sin~\alpha_g sin~\beta_g & cos~\gamma_g 
	\end{bmatrix}
}
\end{equation}
If symmetry rotation matrices $R_g$ are given for all $g\in G$, 
one set of Euler angles $(\alpha_g,\beta_g,\gamma_g)$ can be computed, in particular,  
{\scriptsize
	\begin{eqnarray}\label{eq:EulerAngles}
	\beta_g&=&cos^{-1} R_g (3,3),\nonumber\\
	\alpha_g=\begin{cases}
	-cos^{-1} \frac{R_g(2,1)}{R_g(1,1)} &: \beta_g=0\\
	0 &: \beta_g=\pi\\
	tan^{-1} \frac{R_g(3,2)}{R_g(3,1)} &: o.w.
	\end{cases},&~\text{and }&
	\gamma_g=\begin{cases}
	0 &: \beta_g=0\\
	-tan^{-1} \frac{R_g(2,1)}{R_g(1,1)} &: \beta_g=\pi\\
	tan^{-1} \frac{R_g(1,3)}{R_g(2,3)} &: o.w.
	\end{cases}.
	\end{eqnarray}
}

Symmetry rotation matrices for Tetrahedral group and Octahedral group are available on the Bilbao Crystallographic Server\cite{Aroyo2006b}, though permutations are needed in order that the two sets of symmetry rotation matrices share the same multiplication table with their respective irrep matrices provided by Bilbao Crystallographic Server. Note that the permuted symmetry rotation matrices are exactly the 4th irrep matrices for Tetrahedral group, and they are exactly the 5th irrep matrices for Octahedral group. 

For the icosahedral group, Zheng and Doerschuk~\cite{ZhengDoerschukComputersPhysics1995} give a 3-D representation in the form of real orthonormal matrices with determinant $+1$. Let ${R_I}$ represent the Zheng-Doerschuk 3-D matrices. This 3-D rep is exactly the symmetry rotation matrices for an icosahedron that is positioned in a standard orientation (the $z$ axis passes through two opposite vertices and the $xz$ plane includes one edge of the icosahedron). Note that the 2nd and 3rd irreps provided by Liu-Ping-Chen are also 3-D reps, but compared to Zheng-Doerschuk rep, 
each of these two reps has rotation operators listed in a different order. Hence, to apply the corresponding Euler angles in Eq.~\ref{eq:calDdefn}, it is imperative to permute the Zheng-Doerschuk rep to match the rotation operators for either 2nd or 3rd irrep listed by Liu-Ping-Chen. The permuted Zheng-Doerschuk rep, denoted by ${R_I^{s}}$, is similar to Liu-Ping-Chen 2nd or 3rd irrep. In particular, \begin{equation}
\Gamma_{c;I}^{p}=	{V^p}^{-1}{R_I^s}V^p \text{~for $p=2$ or $3$,}
\end{equation} holds, where $\Gamma_{c;I}^{p}$ denotes the $p$th irrep matrices for Icosahedral group, and $V^p$ is the unitary matrix. Using a computer program, the solutions for ${R_I^{p}}$ and ${V^{p}}$ were found for $p=2$ and $3$, which are tabulated in Table 1 and Table 2 (Appendix B), respectively. 

By applying these symmetry rotation matrices and computing Euler angles (Eq.~\ref{eq:EulerAngles}), Wigner $D$ coefficients are determined by its definition in~\cite[p. 92]{rosvall2008maps}.

\subsection{Construction of real orthonormal basis functions}
By applying the projection operator $\calP_{j,k}^p$ to a real spherical harmonic $Z_{l,m}$ with the rotational property in Eq.~\ref{eq:rotateYlmusingWignerD} and using the real irrep matrices $\Gamma_r^p$, Eq.~\ref{eq:GroupProjectorOperation} becomes
\begin{equation}
\calP_{j,k}^p
Z_{l,m}(\frac{\vx}{x})	
=
\frac{d_p}{\Ngroup}
\sum_{g\in G}
\Gamma_r^p(g)_{j,k}
P(g)Z_{l,m}(\frac{\vx}{x})		
=
\sum_{\m=-l}^{+l}
\calD_{j,k,l,m,\m}^p(g)
Y_{l,\m}(\frac{\vx}{x})
\label{eq:usingallYlm}
\end{equation} 
where
\begin{eqnarray} \label{eq:calDdefn}
& &\calD_{j,k,l,m,\m}^p
=
\frac{d_p}{\Ngroup}
\sum_{g\in G}
\Gamma_r^p(g)_{j,k}
M_{\m,m}^l(g)\\
&=&
\left\{
\begin{array}{lr}
\frac{i d_p}{\sqrt{2}\Ngroup}
\sum_{g\in G}
\Gamma_r^p(g)_{j,k}\left(D_{\m,m}^l(g)-(-1)^{-m} D_{\m,-m}^l(g) \right) & : m < 0\\
\frac{d_p}{\Ngroup}
\sum_{g\in G}
\Gamma_r^p(g)_{j,k}D_{\m,0}^l(g) & : m = 0\\
\frac{d_p}{\sqrt{2}\Ngroup}
\sum_{g\in G}
\Gamma_r^p(g)_{j,k}\left(D_{\m,-m}^l(g)+(-1)^m D_{\m,m}^l(g)\right) & : m > 0
\end{array}
\right.\nonumber
.
\end{eqnarray}

For each $(l,~m,~p)$ triple, 
a normalized real basis function $I_{p;l,m}$ ($\zeta$ in Eq.~\ref{eq:BasisDefn} becomes a shorthand for $(l,m)$) of the polyhedral group $G$, computed by Eqs.~\ref{eq:usingallYlm}--\ref{eq:calDdefn}, has the following vector form,
{
	\begin{equation}\label{eq:VectorBasisFunction}
	\basisplm(\frac{\vx}{x})
	=
	\begin{bmatrix}
	\basisplm(\frac{\vx}{x})[1]\\
	\vdots\\
	\basisplm(\frac{\vx}{x})[d_p]
	\end{bmatrix}
	=
	\frac{1}{c^p_{k,l,m}}
	\begin{bmatrix}
	\calP_{1,k}^p Z_{l,m}(\frac{\vx}{x})\\
	\vdots\\
	\calP_{d_p,k}^p Z_{l,m}(\frac{\vx}{x})
	\end{bmatrix}
	=
	\frac{1}{c^p_{k,l,m}}\sum_{\m=-l}^{l}\boldsymbol{\calD}_{l,m,\m}^pY_{l,\m}(\frac{\vx}{x}), 
	\end{equation}
} where $\boldsymbol{\calD}_{l,m,\m}^p=\begin{bmatrix}
\calD_{1,k,l,m,\m}^p\\
\vdots\\
\calD_{d_p,k,l,m,\m}^p
\end{bmatrix}$ and $c^p_{k,l,m}=\sqrt{\sum_{\m=-l}^{l}|{\calD}_{k,k,l,m,\m}^p|^2}\in\mathbb{R}^+$ for some $k\in\{1,\dots,d_p\}$. 
\par
Because $\Gamma^p$, $P(g)$ and $\psi(\vx)$ in Eq.~\ref{eq:GroupProjectorOperation} are all real, the basis functions $I_{p;l,n}$, which is a linear combination of all real parameters, is real. From Eq.~\ref{eq:VectorBasisFunction}, the $d_p$-dimensional real basis function is also a linear combination of spherical harmonics with the weights $\boldsymbol{\calD}_{k,l,m,\m}^p$, so the imaginary parts of weights and of spherical harmonics are balanced out in the summation. 
\par
Let $\boldsymbol{\calD}_{l,m,\m}^p$ be the $\m$th ($\m=-l,\dots,l$) column of a $(2l+1)\times (2l+1)$ coefficient matrix, denoted by $\boldsymbol{\calH}_{l,m}^p\in\mathbb{C}^{(2l+1)\times (2l+1)}$, of which computation is described in Algorithm~1. Then, Eq.~\ref{eq:VectorBasisFunction} has the matrix form $\basisplm(\frac{\vx}{x})=\boldsymbol{\calH}_{l,m}^p Y^l(\frac{\vx}{x})$, for $\forall \vx\in\mathbb{R}^3$. Note that Algorithm~1 computes $2l+1$ coefficient matrices by varying $m$ through the set $\{-l,\dots,+l\}$, so that $2l+1$ $d_p$-dimensional vectors of basis functions, which are more than necessary for a basis, are computed. Through Gram-Schmidt orthogonalization, the set of coefficient matrices, $\boldsymbol{\calH}_{l,m}^p$ for $m=-l,\dots,l$, shrinks to a smaller set of coefficient matrices, $\boldsymbol{\calH}_{l,n}^p$ for $n=1,\dots,\Npl<2l+1$, and thereby a $\Npl$ dimensional real orthonormal basis is formed. In particular, the final real orthonormal basis functions, denoted by $\basispln(\frac{\vx}{x})$, can be computed by \begin{equation}\label{eq:MatrixBasisFunction}
\basispln(\frac{\vx}{x})
=
\boldsymbol{\calH}_{l,n}^p Y^l(\frac{\vx}{x}),\text{~for $n=1,\dots,\Npl$.}
\end{equation} 
The software of computing real basis functions for polyhedral groups is available from the author. 
\begin{figure}[h]
	\begin{center}
			\includegraphics[width=12.5cm]{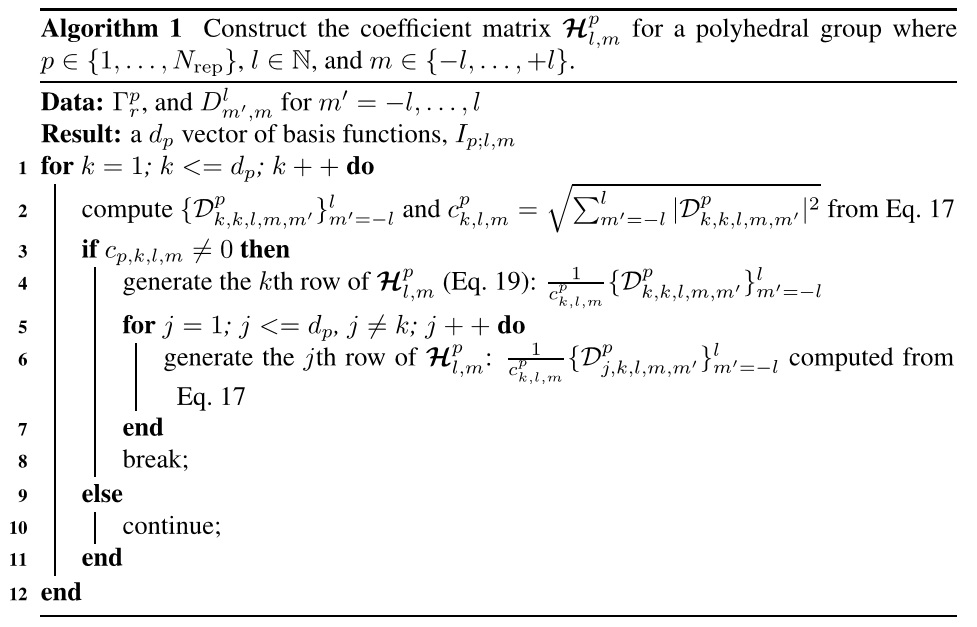}
	\end{center}
	\label{alg:constructbasis}
\end{figure}
\subsection{Properties of the real orthonormal basis functions}\label{sec:RealProperty}
Basis functions of a group, computed by the framework of Eq.~\ref{eq:GroupProjectorOperation}, are guaranteed by Theorem II in~\cite[p.93]{Cornwell1984} to have
\begin{eqnarray}
(\calP_{j,k}^p Z_{l,n},\calP_{j',k'}^{p'} Z_{l',n'})&=&\int \calP_{j,k}^p  \calP_{j',k'}^{p'}Z_{l,n} Z_{l',n'} d\vx=\delta_{pp'}\delta_{kj'}(\calP_{j,k'}^{p}Z_{l,n},Z_{l',n'}),\\
\text{and~}	\calP_{n,n}^p\psi(\vx)&=&\sum_{\zeta}a_{p;\zeta} I_{p;\zeta}
\end{eqnarray}
Therefore, basis functions $I_{p;l,n}$ are orthogonal to each other between different $p$ reps. Also, the subspace spanned by all basis functions $I_{p;l,n}$ for the $p$th rep, denoted by $L_2^p$, are disjoint subspaces with their union to be the $L_2$ space.
In other words, the complete set of basis functions, $I_{p;l,n}$ for $n=1,\dots,\Npl$, $p=1,\dots,\Nirrep$, and $l=0,1,\dots$, form a complete orthonormal system (C.O.N.S.) on $L_2$. Similar results are deducted by~\cite{zheng2000explicit}. Because real basis functions of Tetrahedral group only exist for the 1st and the 4th reps, so these real basis functions only span the corresponding $L_2^1$ and $L_2^4$ subspaces. We then have the following corollary:
\begin{quote}
	Corollary 1: {\sl There exists real basis functions of Octahedral or Icosahedral group which form a C.O.N.S. on $L_2$. There exists no such real basis functions of Tetrahedral group which form a C.O.N.S. on $L_2$.}
\end{quote}

In summary, the real orthonormal basis functions $I_{p;l,n}$, have the following properties:
\begin{enumerate}
	\item
	Each $I_{p;l,n}$ is a $d_p$-dimensional real-valued vector function, i.e.,
	$I_{p;l,n}\in \mathbb{R}^{d_p}$.
	\item
	Each $I_{p;l,n}$ function has a specific transformation property under rotations from the group, defined by Eq.~\ref{eq:BasisDefn}.	
	\item
	The $I_{p;l,n}$ functions are orthonormal on the surface of the sphere.
\end{enumerate} 
The additional properties of real basis functions for Octahedral and Icosahedral groups are
\begin{enumerate}
	\item[4.]
	The subspace of square integrable functions on the surface of the sphere defined by spherical harmonics of index $l$, contains a set of
	$I_{p;l,n}$ functions with a total of $2l+1$ components.
	\item[5.]
	The family of $I_{p;l,n}$ is a complete basis for square integrable functions on the surface of the sphere.
\end{enumerate}
\section{Numerical Results}\label{sec:results}
For Tetrahedral group, the coefficient matrices $\boldsymbol{\calH}_{l,m}^p$ for degree $l=1,\dots,45$, $p=1,4$ and $n=1,..., \Npl$, were computed. The total number of rows of coefficient matrices is $N_{p=1;l}+N_{p=4;l}<2l+1$ for each $l$, which is in agreement with the incompleteness of the orthonormal system formed by Tetrahedral basis on real space (Corollary 1). The basis functions computed by Eq.~\ref{eq:MatrixBasisFunction} are real and orthonormal, which have dimension 1 for $p=1$ and dimension 3 for $p=4$, respectively. 
\par
For both Octahedral and Icosahedral group, the coefficient matrices $\boldsymbol{\calH}_{l,n}^p$ for degree $l=1,\dots,45$, $p=1,\dots,N_{rep}$ and $n=1,..., \Npl$, were computed. For each $l$, by concatenating the rows of all $\sum_{p=1}^{N_{rep}}\Npl$ coefficient matrices into one matrix, a full coefficient matrix $\mathbf{H}^{l}$ was formed. Matrix $\mathbf{H}^{l}$ for all $l$'s were verified to be unitary and to have the dimension $(2l+1)\times (2l+1)$, i.e., $\sum_{p=1}^{N_{rep}}d_p*\Npl=2l+1$, so that the set of real basis functions computed by Eq.~\ref{eq:MatrixBasisFunction} that spans the subspace of square integrable functions on the surface of the sphere defined by degree $l$ form an (2l+1) orthonormal basis. 
\par
As a result, properties described in Section~\ref{sec:RealProperty} had been verified for the corresponding polyhedral group.

\begin{figure}[h]
	\begin{center}
		\begin{tabular}{cc}					
			\includegraphics[width=2cm]{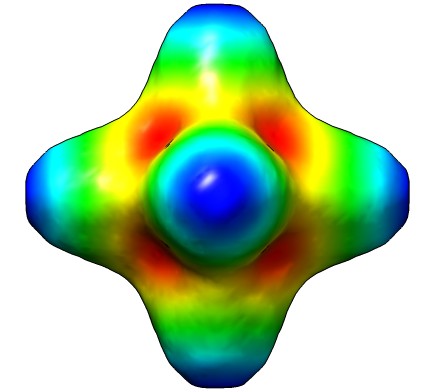}
			&
			\includegraphics[width=2cm]{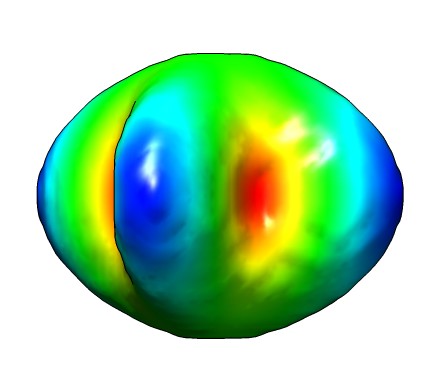}
			\\
			$T_{1;4,1,1}$			
			&
			$T_{4;4,1,1}$
		\end{tabular}	
	\end{center}	
	\caption{
		\label{fig:irredrepT}
		Examples of the real basis functions of the Tetrahedral group. }
\end{figure}

\begin{figure}[h]
	\begin{center}
		\begin{tabular}{ccccc}
			\includegraphics[width=1.8cm]{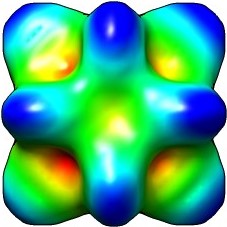}
			&
			\includegraphics[width=2cm]{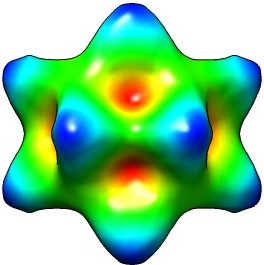}
			&
			\includegraphics[width=2cm]{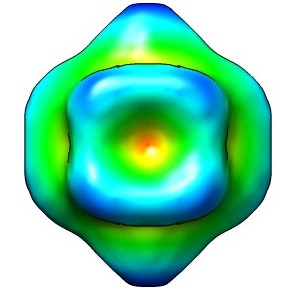}
			&
			\includegraphics[width=2cm]{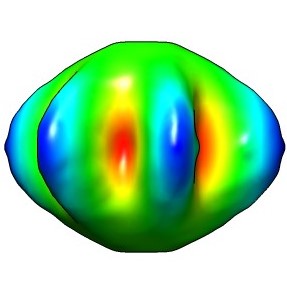}
			&
			\includegraphics[width=2cm]{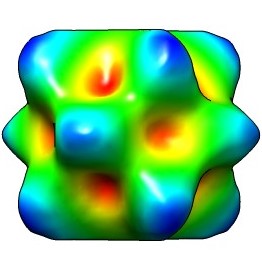}
			\\
			$O_{1;6,1,1}$
			&
			$O_{2;6,1,1}$
			&
			$O_{3;6,1,2}$
			&
			$O_{4;6,1,1}$
			&
			$O_{5;6,1,1}$
		\end{tabular}		
	\end{center}
	\caption{\label{fig:irredrepO}
		Examples of the real basis functions of the Octahedral group, denoted by $O_{p,l,n,j}$.  
	}
\end{figure}

\begin{figure}[h!]
	\begin{center}
		\begin{tabular}{ccccc}
			\includegraphics[width=2cm]{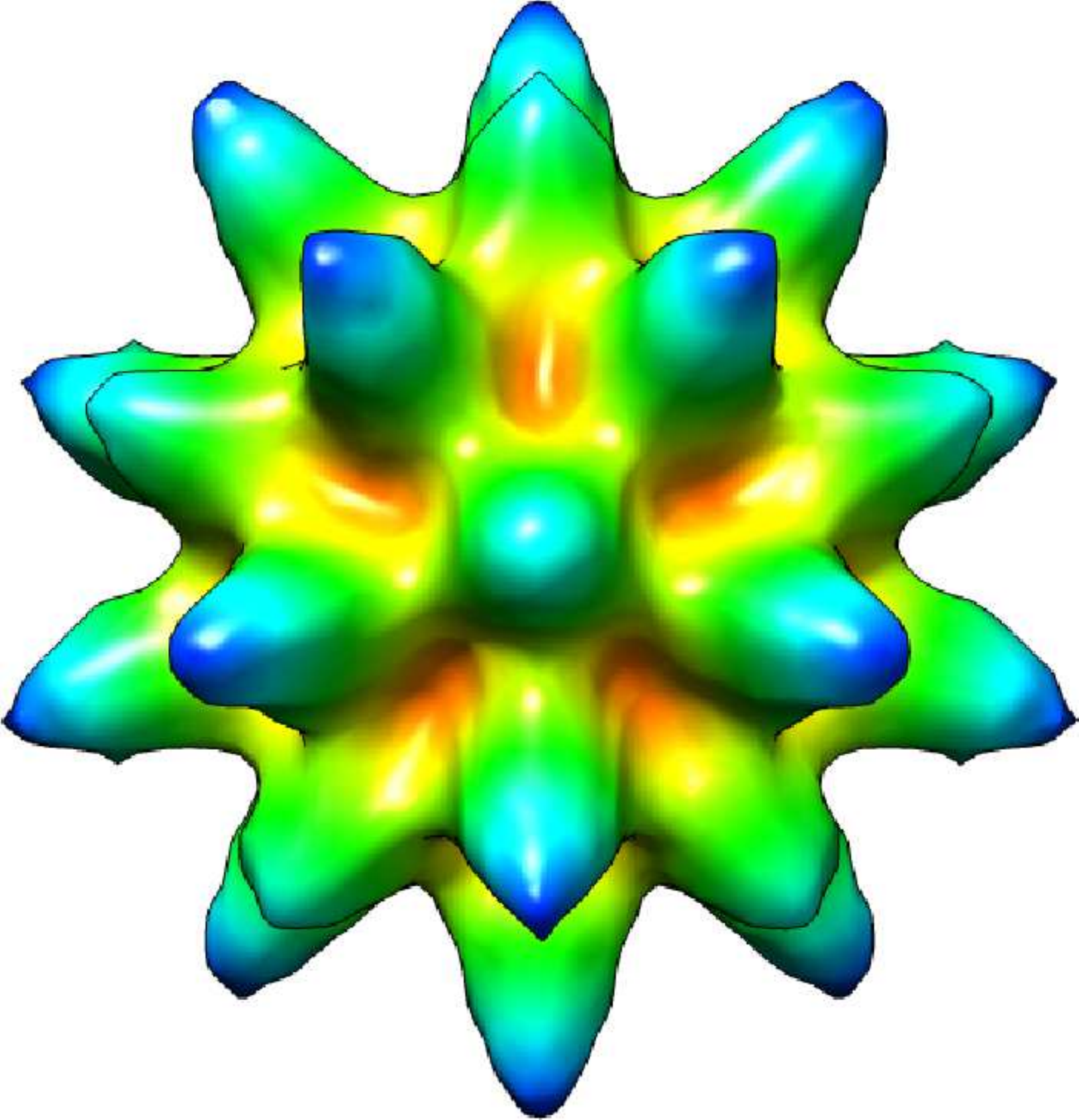}
			&
			\includegraphics[width=2cm]{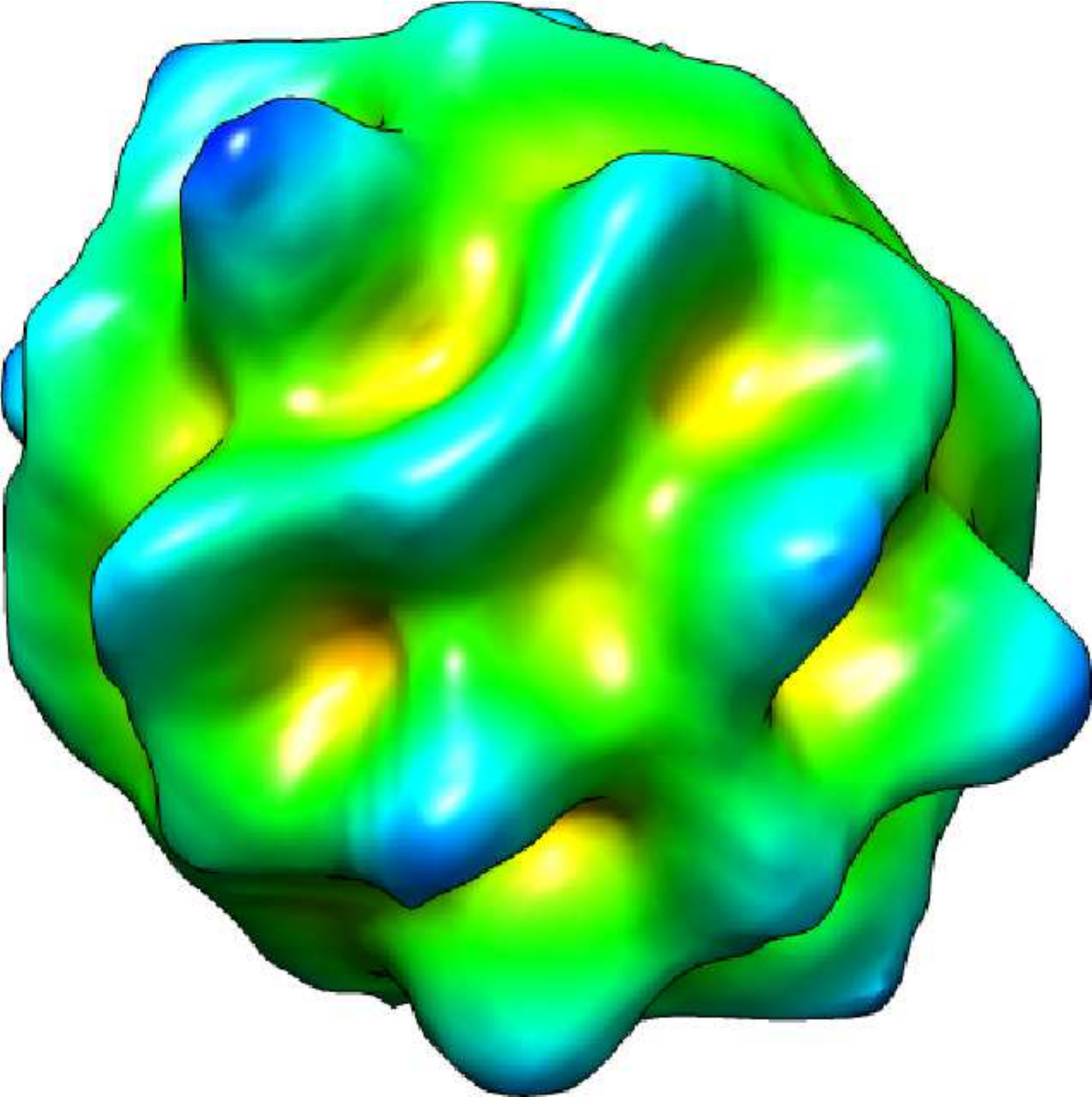}
			&
			\includegraphics[width=2cm]{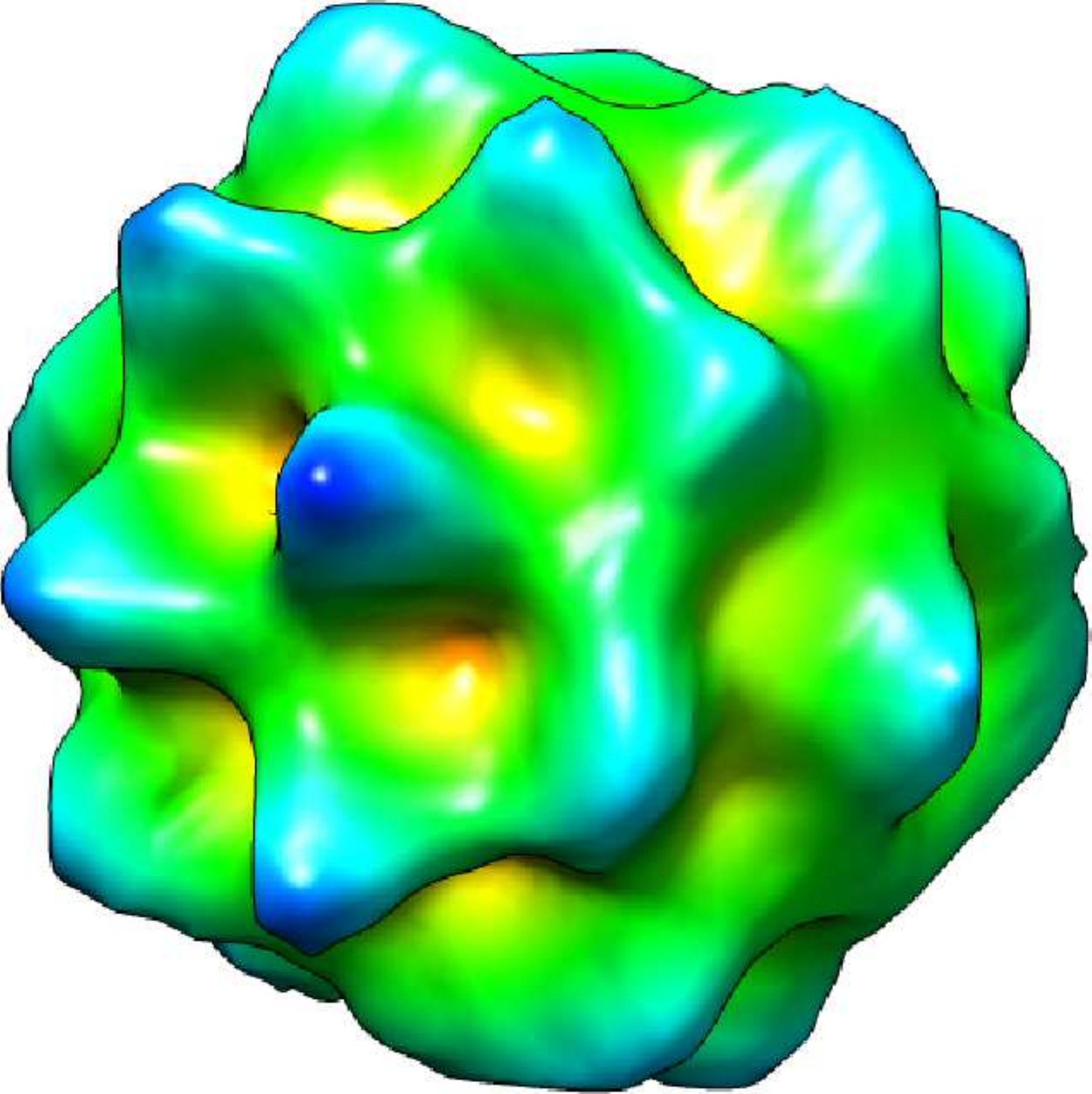}
			&
			\includegraphics[width=2cm]{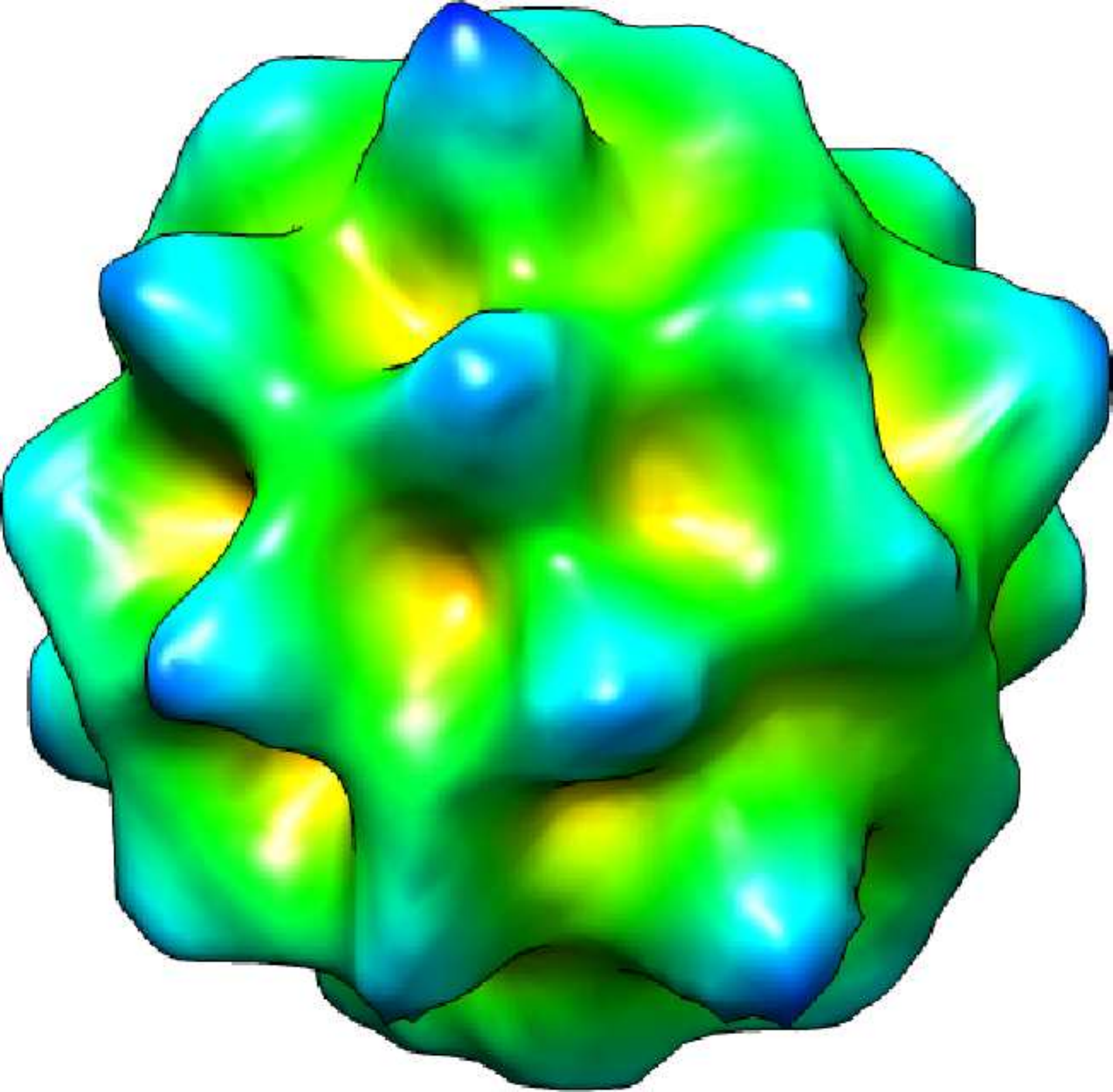}
			&
			\includegraphics[width=2cm]{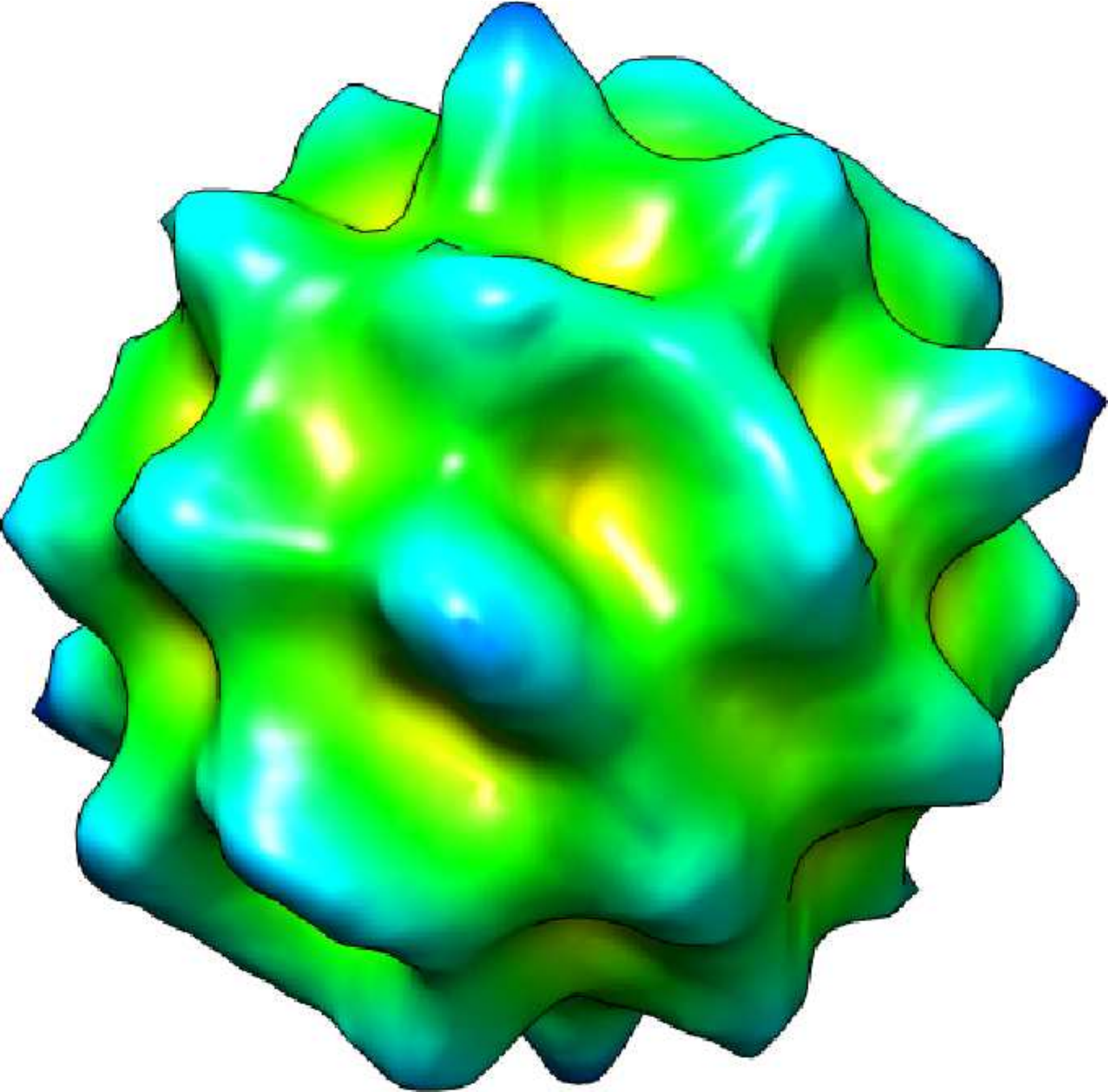}
			\\
			$I_{1;10,1,1}$
			&
			$I_{2;10,1,1}$
			&
			$I_{3;10,1,1}$
			&
			$I_{4;10,1,1}$
			&
			$I_{5;10,1,1}$
		\end{tabular}	
	\end{center}
	\caption{\label{fig:irredrepI}
		Examples of the real basis functions of the icosahedral group, denoted by $I_{p,l,n,j}$. 
	}
\end{figure}

For visualization purposes, define the 3-D object $\xi(\vx)$ by 
\begin{equation}
\xi(\vx)
=
\left\{
\begin{array}{ll}
1 , & x\le \kappa_1+\kappa_2 \scalarbasispln(\vx/x) \\
0 , & \mbox{otherwise}
\end{array}
\right.
\label{eq:visualizationofharmonic}
\end{equation}
where $\kappa_1$ and $\kappa_2$ are chosen so that $\kappa_1+\kappa_2
\scalarbasispln(\vx/x)$ varies between 0.5 and 1. Examples of real basis functions of Tetrahedral, Octahedral and Icosahedral groups for $l$ = 15 are shown in Figure~\ref{fig:irredrepT}, Figure~\ref{fig:irredrepO} and ~\ref{fig:irredrepI}, respectively.  The surfaces of 3-D objects defined by
Eq.~\ref{eq:visualizationofharmonic} are visualized by UCSF
Chimera~\cite{PettersenHuangCouchGreenblattMengFerrin2004}
where the color indicates the distance from the center of
the object. Note that the p = 1 exhibits all of the symmetries of the solid.

\section{Aknowledgement}
I am grateful to Prof. Peter C. Doerschuk who support me to do this work. Peter is the person who introduced me to my first group theory book, "Group Theory in Physics" written by Cornwell. Before I started my work, Peter had the idea of generating complex basis functions of the icosahedral group through the standard group projection approach. He brought Liu-Ping-Chen's paper to me to start with. Appendix C and the visualization of real basis functions (Eq. 22) are Peter's idea and contribution. 

I am especially grateful to Professor Dan Mihai Barbasch for his willingness of discussing my group theory problem with me and leading me to the right direction of solving this representation theory problem. Dan was the only professor in Math Department who took a serious thinking on my problem when I stop by Cornell Math Department in that week. Moreover, he sent me several hand written notes about his thoughts in the following couple weeks. The method of generating real irrep matrices in Section 2.2 is purely from the note written by Dan. 

I am also grateful to my friends Hung Tran and Chenxi Wu from Department of Mathematics for their helpful discussions on group theory. Hung is the person who led me to talk to everyone in Math Department about my problem until we finally reached Dan. I would like to thank my friend Chaoxu Tong from ORIE for his helpful discussion on linear algebra, especially finding basis for a random matrix. I would also like to thank my best friend Yuguang Gao for all her support and her critiques on writing the applied math paper. 

Finally, I am full of gratitude to my husband, Dr. Andy Sun, for his critiques on an earlier version of this paper, which motivated me to explore further about the relationship between the real irrep matrices and real basis functions, which brought the current Proposition 1 and 2. My discussion with him helped me make the whole story more complete with no logic holes. 

\clearpage
\section{Appedix A}
One approach to computing $S^p$ is described in this section. First, we introduce Lemma 1, which shows the relationship between the coneigenvectors of a random non-singular complex symmetric matrix and the eigenvectors of its real representation matrix.

Lemma 1: Let $A\in\mathbb{C}^{n\times n}$ be a non-singular symmetric matrix with $A\overline{A}=I$. Suppose $A=A_1+iA_2$, where $A_1$, $A_2\in\mathbb{R}^{n\times n}$ is the real and imaginary part of matrix $A$, receptively. Consider the matrix $A$'s real structure, $B=\bigl[\begin{smallmatrix}
A_1 & A_2 \\ -A_2 & A_1 
\end{smallmatrix} \bigr]\in\mathbb{R}^{2n\times 2n}$. Then, the following three properties hold:

1) $B$ is nonsingular.

2)$B \begin{bmatrix}
x \\ 
-y 
\end{bmatrix}=\lambda \begin{bmatrix}
x \\ 
-y 
\end{bmatrix}
$ if and only if $B\begin{bmatrix}	y \\ 
x
\end{bmatrix}=-\lambda \begin{bmatrix}
y \\ 
x 
\end{bmatrix}$ so the eigenvalues of B  appear in $(+, -)$  pairs.

3)Let $\begin{bmatrix} x_1 \\ -y_1 \end{bmatrix},\dots, \begin{bmatrix} x_n \\ -y_n \end{bmatrix}$ be the orthonormal eigenvectors of $B$ associated with its positive eigenvalues $\lambda_1,\dots,\lambda_n$.  Then, $\{x_1-iy_1, \dots, x_n-iy_n\}$ is the set of orthonormal coneigenvectors of $A$ associated with coneigenvlues $\lambda_1,\dots,\lambda_n$. In other words, $A\overline{(x_k-iy_k)}=\lambda_k(x_k-iy_k)$ for $k=1,...,n$.

Proof: It's straightforward that matrices $A_1$, $A_2$ and $B$ are all real and symmetric. Hence, eigenvectors of $B$ (if existing) are real as well. 

1)Note that $M=\begin{bmatrix}I&-iI\\ 0 & I\end{bmatrix} \begin{bmatrix}A_1 & A_2\\ A_2& -A_1\end{bmatrix} \begin{bmatrix}I & 0\\ iI & I\end{bmatrix}= \begin{bmatrix}0 & A_2+iA_1\\ A_2-iA_1 & -A_1\end{bmatrix}$ and $\det(B)=\det(M)=\det((A_2+iA_1)(A_2-iA_1)-0(-A_1))=\det(A\overline{A})=
|\det(A)|^2>0$ for nonsingular $A$, which implies $B$ is nonsingular and has $2n$ distinct eigenvalues. 

2) The left equation holds $\Leftrightarrow\left\{\begin{array}{l}
A_1 x-A_2 y=\lambda x \\
A_2 x+A_1 y=-\lambda y
\end{array}\right.\Leftrightarrow\left\{\begin{array}{l}
A_2 y-A_1 x=-\lambda x \\
A_2 x+A_1 y=-\lambda y
\end{array}\right.\Leftrightarrow$ the right equation holds.

3) Define matrices $X=\begin{bmatrix} x_1 & \dots & x_n \end{bmatrix}\in\mathbb{R}^{n\times n}$, $Y=\begin{bmatrix} y_1 & \dots & y_n \end{bmatrix}\in\mathbb{R}^{n\times n}$, $\Sigma=diag(\lambda_1,....\dots,\lambda_n)\in\mathbb{R}^{n\times n}$,
and let $U=X - iY$. Then, we have matrix form equation $B\begin{bmatrix} X \\ -Y \end{bmatrix}=\begin{bmatrix} A_1 & A_2 \\ A_2 & -A_1 \end{bmatrix}\begin{bmatrix} X \\ -Y \end{bmatrix}=\begin{bmatrix} X \\ -Y \end{bmatrix}\Sigma$ which is equivalent to $\left\{\begin{array}{l}
A_1 X-A_2 Y=X\Sigma \\
A_2 X+A_1 Y=-Y\Sigma 
\end{array}\right.$. Multiplying the second equation by $i$ and then adding to the first equation gives $(X-iY)\Sigma=(A_1X-A_2Y)+i(A_2X+A_1Y)=(A_1+iA_2)X+(iA_1-A_2)Y=(A_1+iA_2)X+(A_1+iA_2)iY=(A_1+iA_2)(X+iY)$. That is $A\overline{U}=U\Sigma$. Thus, $U=X-iY$ is an unitary matrix columned by coneigenvectors of $A$ and diagonal entries of $\Sigma$ are coneigenvalues of $A$. The unitarity of $U$ due to the orthonormal property of eigenvectors of $B$. Thus, $A=U\Sigma\overline{U}^{-1}=U\Sigma{U}^{T}$.  
$\qedwhite$

Followed by Lemma 1, the orthonormal eigenvectors of the real representation matrix of $C^p$ and thereby the coneigenvectors of $C^p$ are computed. Consider the matrix $\begin{bmatrix} X^p \\ -Y^p \end{bmatrix}$ columned by the orthonormal eigenvectors of the real representation matrix of $C^p$ that associates to eigenvalue matrix $\Sigma^p$. Let $S^p=X^p-iY^p$. Then, $C^p\overline{S^p}=S^p\Sigma^p$.
By Proposition 4.6.6 in~\cite[p.~246]{HornJohnson1985}, $\Sigma^p=I$ since $C^p \overline{C^p}=I$. Thus, $C^p\overline{S^p}=S^p$, or equivalently, $C^p=S^p{S^p}^T$. 

With this $S^p$, the LHS of Eq. \ref{eq:realProperty} is equivalent to ${S^p}^T\overline{\Gamma_c^p}\overline{S^p}={S^p}^T(\overline{C^p}\Gamma_c^p C^p)\overline{S^p}={S^p}^T\overline{C^p}\Gamma_c^p (C^p\overline{S^p})={S^p}^T\overline{{C^p}^T}\Gamma_c^pS^p=(\overline{C^p}{S^p})^T\Gamma_c^pS^p
=({C^p}\overline{S^p})^H\Gamma_c^pS^p={S^p}^H\Gamma_c^pS^p$, which is the RHS. The first equivalent sign is based on Eq. \ref{eq:complexSymSimilarity}; the third equivalence comes from the symmetry of $C^p$; the fifth equal sign due to the fact that both $C^p$ and $S^p$ are unitary.

\section{Appendix B}
\begin{quote}
	Proposition 1: {\sl Real orthonormal basis functions, i.e., $I_{p;\zeta}\in\mathbb{R}^{d_p}$, generate real orthogonal irrep matrices of a group. }
\end{quote} 
Proof: Given any $\zeta\in\{1,\dots,N_\zeta\}$, let $I_{p;\zeta}[i]$ and $I_{p;\zeta}[j]$ be the $i$th and the $j$th entry of the basis function vector $I_{p;\zeta}$ for $i,j\in\{1,\dots,d_p\}$. Following Eq.(5-73) in~\cite[p.91]{bishop2012group}, 
\[(P(g)I_{p;\zeta}[i],P(g)I_{p;\zeta}[j])=(I_{p;\zeta}[i],I_{p;\zeta}[j])\] For basis functions that are real and orthonormal, we have $\overline{I_{p;\zeta}}=I_{p;\zeta}$, $\overline{P(g)I_{p;\zeta}}=P(g)I_{p;\zeta}$ and $(I_{p;\zeta}[i],I_{p;\zeta}[j])=\delta_{i,j}$. Then, for any $g\in G$,
\begin{eqnarray}\label{eq:LHSpfProp1}
LHS&=&\int\overline{P(g)I_{p;\zeta}[i]} 
P(g)I_{p;\zeta}[j]d\vx=\int{P(g)I_{p;\zeta}[i]} 
P(g)I_{p;\zeta}[j]d\vx\nonumber\\
&=&\int{\sum_{h=1}^{d_p}\Gamma^p(g)_{h,i}I_{p;\zeta}[h] \sum_{k=1}^{d_p}\Gamma^p(g)_{k,j}I_{p;\zeta}[k]}d\vx\nonumber\\
&=&\sum_{h=1}^{d_p}\sum_{k=1}^{d_p}\Gamma^p(g)_{h,i}\Gamma^p(g)_{k,j}\int{I_{p;\zeta}[h] I_{p;\zeta}[k]}d\vx\nonumber\\
&=&\sum_{h=1}^{d_p}\sum_{k=1}^{d_p}\Gamma^p(g)_{h,i}\Gamma^p(g)_{k,j}\int{\overline{I_{p;\zeta}[h]} I_{p;\zeta}[k]}d\vx\nonumber\\
&=&\sum_{h=1}^{d_p}\sum_{k=1}^{d_p}\Gamma^p(g)_{h,i}\Gamma^p(g)_{k,j}(I_{p;\zeta}[h],I_{p;\zeta}[k])=\sum_{h=1}^{d_p}\sum_{k=1}^{d_p}\Gamma^p(g)_{h,i}\Gamma^p(g)_{k,j}\delta_{h,k}\nonumber\\
&=&\sum_{k=1}^{d_p}\Gamma^p(g)_{k,i}\Gamma^p(g)_{k,j}=RHS=\delta_{i,j}.
\end{eqnarray}
The matrix form of the Eq.~\ref{eq:LHSpfProp1} gives $\Gamma^p(g)^T\Gamma^p(g)=I$. Therefore $\Gamma^p(g)$ is a real orthogonal matrix.

\section{Appendix C}
\begin{center}
	\begin{tabular}{rrrrrrrrrr}
		1 & 2 & 3 & 4 & 5 & 6 & 7 & 8 & 9 & 10 \\
		1 & 2 & 5 & 9 & 17 & 10 & 27 & 13 & 21 & 18 \\
		\hline
		11 & 12 & 13 & 14 & 15 & 16 & 17 & 18 & 19 & 20 \\
		24 & 15 & 26 & 3 & 4 & 48 & 45 & 56 & 54 & 49 \\
		\hline
		21 & 22 & 23 & 24 & 25 & 26 & 27 & 28 & 29 & 30 \\
		60 & 36 & 52 & 42 & 38 & 14 & 16 & 47 & 40 & 46 \\
		\hline
		31 & 32 & 33 & 34 & 35 & 36 & 37 & 38 & 39 & 40 \\
		55 & 41 & 53 & 20 & 29 & 6 & 12 & 57 & 39 & 8 \\
		\hline
		41 & 42 & 43 & 44 & 45 & 46 & 47 & 48 & 49 & 50 \\
		22 & 44 & 58 & 28 & 25 & 11 & 31 & 59 & 33 & 30 \\
		\hline
		51 & 52 & 53 & 54 & 55 & 56 & 57 & 58 & 59 & 60 \\
		19 & 43 & 35 & 34 & 37 & 23 & 7 & 50 & 32 & 51
	\end{tabular}
\end{center}
\begin{equation}
V^{p=2}
=
\left[
\begin{array}{cccc}
-1/\sqrt{2} & 0 & -1/\sqrt{2} \\
-i/\sqrt{2} & 0 & i/\sqrt{2} \\
0 & 1 & 0
\end{array}
\right]
.
\end{equation}

\begin{center}
	\begin{tabular}{rrrrrrrrrr}
		1 & 2 & 3 & 4 & 5 & 6 & 7 & 8 & 9 & 10 \\
		1 & 4 & 3 & 36 & 52 & 38 & 42 & 49 & 60 & 54 \\
		\hline
		11 & 12 & 13 & 14 & 15 & 16 & 17 & 18 & 19 & 20 \\
		56 & 48 & 45 & 2 & 5 & 24 & 18 & 15 & 26 & 21 \\
		\hline
		21 & 22 & 23 & 24 & 25 & 26 & 27 & 28 & 29 & 30 \\
		13 & 10 & 27 & 17 & 9 & 46 & 55 & 22 & 8 & 25 \\
		\hline
		31 & 32 & 33 & 34 & 35 & 36 & 37 & 38 & 39 & 40 \\
		28 & 20 & 29 & 53 & 41 & 40 & 47 & 12 & 6 & 39 \\
		\hline
		41 & 42 & 43 & 44 & 45 & 46 & 47 & 48 & 49 & 50 \\
		57 & 16 & 14 & 44 & 58 & 50 & 31 & 11 & 32 & 43 \\
		\hline
		51 & 52 & 53 & 54 & 55 & 56 & 57 & 58 & 59 & 60 \\
		51 & 19 & 33 & 35 & 7 & 59 & 37 & 23 & 34 & 30
	\end{tabular}
\end{center}
\begin{equation}
V^{p=3}
=
\left[
\begin{array}{cccc}
-1/\sqrt{2} & 0 & -1/\sqrt{2} \\
i/\sqrt{2} & 0 & -i/\sqrt{2} \\
0 & 1 & 0
\end{array}
\right]
.
\end{equation}

\clearpage

\bibliographystyle{plain}
\bibliography{referencesnB}

\end{document}